\documentclass[12pt, oneside]{amsart}   	
\usepackage{geometry}                		
\usepackage{graphicx}				
\usepackage{amsmath,amssymb,amsthm,amscd}
\usepackage{ascmac}
\usepackage{slashbox}
\usepackage{indentfirst}
\usepackage[all]{xy}
\usepackage[OT2,T1]{fontenc}
\DeclareSymbolFont{cyrletters}{OT2}{wncyr}{m}{n}
\DeclareMathSymbol{\Sha}{\mathalpha}{cyrletters}{"58}
\DeclareMathOperator{\Gal}{\rm{Gal}}
\DeclareMathOperator{\dim_Fp}{{\rm{dim}}_{\mathbb{F}_{\textit{p}}}}
\DeclareMathOperator{\Cl_K}{{\rm{Cl}}_{\textit{K}}}
\DeclareMathOperator{\Sym^j}{{\rm{Sym}}^{\textit{j}}}
\newtheorem{thm}{\textbf{Theorem}}[section]
\newtheorem*{thm*}{\textbf{Theorem}}

\newtheorem{defn}[thm]{\textbf{Definition}}
\newtheorem{prop}[thm]{\textbf{Proposition}}
\newtheorem{lem}[thm]{\textbf{Lemma}}

\theoremstyle{definition}
\newtheorem{rem}[thm]{\textbf{Remark}}
\newtheorem*{rem*}{\textbf{Remark}}

\title[]{Ideal class groups of number fields and Bloch-Kato's Tate-Shafarevich groups for symmetric powers of elliptic curves}
\author{Naoto Dainobu}
\address{Department of Mathematics \\
3-14-1 Hiyoshi, Kohoku-ku, Yokohama-shi, Kanagawa 223-8522 Japan}
\email{vicarious@keio.jp}

\begin{document}
\maketitle
\begin{abstract}
 For an elliptic curve $E$ over $\mathbb{Q}$, putting $K=\mathbb{Q}(E[p])$ which is the $p$-th division field of $E$ for an odd prime $p$, we study the ideal class group $\Cl_K$ of $K$ as a $\Gal(K/\mathbb{Q})$-module. More precisely, for any $j$ with $1\leqslant j \leqslant p-2$, we give a condition that $\Cl_K\otimes \mathbb{F}_p$ has the symmetric power $\Sym^jE[p]$ of $E[p]$ as its quotient $\Gal(K/\mathbb{Q})$-module, in terms of Bloch-Kato's Tate-Shafarevich group of $\Sym^j V_p E$. Here $V_p E$ denotes the rational $p$-adic Tate module of $E$. This is a partial generalization of a result of Prasad and Shekhar for the case $j=1$.
\end{abstract}

\section{Introduction}
 The ideal class groups of number fields, the Tate-Shafarevich groups and the Selmer groups of elliptic curves are central objects to study in number theory. Many people have noticed the existence of various relations between the class groups and the Tate-Shafarevich groups, or the class groups and the Selmer groups. For example in \cite{W}, Washington considered a specific elliptic curve defined by the equation of the simplest cubic, and studied a relation between its 2-Selmer group and the class group of its 2-division field. In \cite{Nek}, Nekov\'a\v{r} studied a relation between the ideal class groups of certain quadratic fields and the Tate-Shafarevich groups of twists of the cubic Fermat curve. We note here that they studied the ideal class groups of abelian number fields over $\mathbb{Q}$. In this paper, for an elliptic curve $E$ over $\mathbb{Q}$ and an odd prime $p$, we suppose that the group of $p$-torsion points $E[p]$ of $E$ is irreducible as a Galois module, and study the ideal class group of the $p$-th division field $K=\mathbb{Q}(E[p])$ of $E$ which is a \textit{non-commutative} Galois extension of $\mathbb{Q}$. More precisely, we relate the ideal class group $\Cl_K$ of $K$ with Bloch-Kato's Tate-Shafarevich groups for symmetric powers of $V_p E$, where $V_p E$ denotes the rational $p$-adic Tate module of $E$.

 Recently Prasad and Shekhar have proved the following theorem on $\Cl_K$ with $K=\mathbb{Q}(E[p])$, which we first recall. In the situation above, the Galois group $G:=\Gal(K/\mathbb{Q})$ acts on the class group $\Cl_K$. In \cite{PS}, they considered $\Cl_K$ as a $G$-module and proved the following result relating $\Cl_K$ with the Tate-Shafarevich group $\Sha(E/\mathbb{Q})$ of $E$ over $\mathbb{Q}$.
\begin{thm*}[Prasad-Shekhar]\label{PS}
Let $\rho_{E, p}:\Gal (\bar{\mathbb{Q}}/\mathbb{Q}) \rightarrow  \mathrm{Aut}(E[p]) \cong \mathrm{GL}_2(\mathbb{F}_p)$ be the $\mathbb{F}_p$-valued Galois representation associated to $E$. Suppose that the following conditions on $E$ hold:
\begin{itemize}
\item[$(a)$] $E$ has good reduction at $p$.
\item[$(b)$] In the case that $E$ has good ordinary reduction at $p$, $a_p (E) \equiv 1 \pmod p$, and $E$ has no CM over an extension of $\mathbb{Q}_p$, then $\rho_{E, p}$ is wildly ramified at $p$.
\item[$(c)$] For every prime number $l \neq p$, the Tamagawa number $c_{l}(E/\mathbb{Q}_{l})$ of $E/\mathbb{Q}_{l}$ is prime to $p$. 
\item[$(d)$] $E[p]$ is an irreducible $\Gal (\bar{\mathbb{Q}}/\mathbb{Q})$-module.
\end{itemize}
Then the condition $\dim_Fp(\Sha(E/\mathbb{Q})[p])\geqslant 2$ implies that the $\mathbb{F}_p$-representation $\Cl_K \otimes \mathbb{F}_p$ of $G$ has $E[p]$ as its quotient representation. 
\end{thm*}

From the above theorem, we see that the $\mathbb{F}_p$-rank of $\Sha(E/\mathbb{Q})[p]$ gives us the information on $\Cl_K \otimes \mathbb{F}_p$ as a $G$-module. We remark that they also studied in \cite{PS} the condition on which $\Cl_K \otimes \mathbb{F}_p$ has $E[p]$ as its quotient representation even if $\Sha(E/\mathbb{Q})[p]=0$.

The main result of this article is an analogy of the above theorem for the symmetric powers of $E[p]$. In the following, we further assume that the representation $\rho_{E, p}$ is surjective, so the Galois group $G=\Gal(K/\mathbb{Q})$ is isomorphic to $\mathrm{GL}_2(\mathbb{F}_p)$. It is a well-known fact that any irreducible representation of $G=\mathrm{GL}_2(\mathbb{F}_p)$ in characteristic $p$ is of the form $\Sym^j E[p]\otimes \mathrm{det}^i$ $(0\leqslant j \leqslant p-1, \ 0\leqslant i \leqslant p-2$), where $\mathrm{det}$ denotes the determinant character of $\mathrm{GL}_2(\mathbb{F}_p)$. So taking the above theorem one step further, we consider the condition on which the $\mathbb{F}_p$-representation $\Cl_K \otimes \mathbb{F}_p$ of $G$ has an irreducible representation $\Sym^j E[p]$ as its quotient representation. 

Now we explain our main result. For any $j$ with $1\leqslant j \leqslant p-1$, we define $V^j_p:=\Sym^j (V_p E)$ to simplify the notation. One of the key objects in the main result is Bloch-Kato's Tate-Shafarevich group $\Sha^{BK}_{\mathbb{Q}}(V^j_p)$ of $V^j_p$ whose definition we shall recall in Definition \ref{Sha} in Section 2. See also \cite[Definition 5.1]{BK}. The main result of this article is as follows.
\begin{thm}\label{main}
Let $p>3$. For any $j$ with $1\leqslant j \leqslant p-2$, suppose that the following conditions on $E$ hold:
\begin{itemize}
\item[$(a^{\prime})$] $E$ has good reduction at $p$.
\item[$(b^{\prime})$] In the case that $E$ has good ordinary reduction at $p$, $a_p (E)^j \equiv 1\pmod p$, and $E$ has no CM over an extension of $\mathbb{Q}_p$, then $\rho_{E, p}$ is wildly ramified at $p$.
\item[$(c^{\prime})$] If $E$ has potentially multiplicative reduction at $l \neq p$, then $v_{l}(j(E))$ is prime to $p$, where $v_l$ denotes the normalized $l$-adic valuation and $j(E)$ the $j$-invariant of $E$.
\item[$(d^{\prime})$] The representation $\rho_{E, p}$ is surjective.
\end{itemize}
Then the condition $\dim_Fp(\Sha^{BK}_{\mathbb{Q}}(V^j_p)[p])\geqslant j+1$ implies that the $\mathbb{F}_p$-representation $\Cl_K \otimes \mathbb{F}_p$ of $G$ has $\Sym^j E[p]$ as its quotient representation.
\end{thm}
\begin{rem}
The assumptions $(a^{\prime})$ in Theorem \ref{main} and $(a)$ in the theorem of Prasad and Shekhar are the same. The assumption $(b^{\prime})$ for $j=1$ is exactly $(b)$. The assumption $(c^{\prime})$ implies the assumption $(c)$ when $E$ has good or split multiplicative reduction at every prime since we have $c_{l}(E/\mathbb{Q}_l)=-v_{l}(j(E))$ when $E$ has split multiplicative reduction at $l$.
The assumption $(d^{\prime})$ also implies the assumption $(d)$. Hence when $j=1$, we can deduce the  theorem of Prasad and Shekhar from our Theorem \ref{main} if $E$ has only good or split multiplicative reduction at any prime and the representation $\rho_{E, p}$ is surjective. So interestingly, using Bloch-Kato's Tate-Shafarevich groups $\Sha^{BK}_{\mathbb{Q}}(V^j_p)$ for various $j$ other than $j=1$, we can get more information about the structure of $\Cl_K \otimes \mathbb{F}_p$ as a $G$-module. 
\end{rem}

\begin{rem}
We can generalize the argument in this paper for more general Galois representations other than $\Sym^j E[p]$. In Remark \ref{rem}, we explain that we can show a result analogous to Theorem \ref{main} for the $\mathbb{F}_p$-valued Galois representations attached to modular forms using a similar argument in this article.
\end{rem}

We give a sketch of the proof of Theorem \ref{main} in Section 2 dividing it into 3 steps. We prove step 1 in Section 3, step 2 in Section 4 and step 3 in Sections 5 and 6.
\section{A sketch of the proof}
We mainly follow the strategy of the proof of Prasad and Shekhar in \cite{PS}. They used the classical $p$-Selmer group $\mathrm{Sel}_{p}(E/\mathbb{Q})$ but, to treat representations such as $\Sym^j E[p]$, we have to deal with Bloch-Kato's Selmer group $H^1_f$ which we first recall.

\vspace{2mm}
For a field $F$, $G_F$ denotes its absolute Galois group $\Gal (\overline{F}/F)$. We define $T^j_p:=\Sym^j(T_p E)$, $A^j_p:=V^j_p/T^j_p\cong \Sym^j E[p^{\infty}]$, where $T_p E$ is the integral $p$-adic Tate module of $E$. For every prime $l$, we define a local condition $H_f^1(\mathbb{Q}_l, V^j_p)$ in $H^1(\mathbb{Q}_l, V^j_p)$ as 
\[
   \begin{cases}
    H_f^1(\mathbb{Q}_l, V^j_p):= \mathrm{Ker}\left(H^1(\mathbb{Q}_l, V^j_p)\rightarrow H^1(\mathbb{Q}^{\mathrm{ur}}_l, V^j_p)\right)   & (l \neq p) \\
    H_f^1(\mathbb{Q}_p, V^j_p):= \mathrm{Ker}\left(H^1(\mathbb{Q}_p, V^j_p)\rightarrow H^1(\mathbb{Q}_p, V^j_p\otimes\mathbf{B}_{\mathrm{crys}})\right) & (l=p).
  \end{cases}
\]
Here $\mathbb{Q}^{\mathrm{ur}}_l$ is the maximal unramified extension of $\mathbb{Q}_l$ and $\mathbf{B}_{\mathrm{crys}}$ denotes Fontaine's crystalline period ring which is defined in \cite[Section 1]{BK}. Then we define $H_f^1(\mathbb{Q}_l, A^j_p):=\pi\left(H_f^1(\mathbb{Q}_l, V^j_p)\right)$ for each prime $l$, where $\pi: H^1 (\mathbb{Q}_l, V^j_p) \rightarrow H^1 (\mathbb{Q}_l, A^j_p)$ is the homomorphism induced by the natural map $\pi: V^j_p \rightarrow A^j_p$. We define Bloch-Kato's Selmer group for $V_p^j$ and $A_p^j$ using these local conditions.
\begin{defn}\label{Sel}
For $V_p^j=\Sym^j (V_p E)$ and $A^j_p=V^j_p/T^j_p\cong \Sym^j E[p^{\infty}]$, we define Bloch-Kato's Selmer groups as
\[
H^1_f(\mathbb{Q}, V^j_p) := \mathrm{Ker}\left(H^1(\mathbb{Q}, V^j_p) \xrightarrow{\prod \mathrm{Loc}_l} \prod_{l} \frac{H^1(\mathbb{Q}_l, V^j_p)}{H^1_f(\mathbb{Q}_l, V^j_p)} \right),
\]
\[
H^1_f(\mathbb{Q}, A^j_p) := \mathrm{Ker}\left(H^1(\mathbb{Q}, A^j_p) \xrightarrow{\prod \mathrm{Loc}_l} \prod_{l} \frac{H^1(\mathbb{Q}_l, A^j_p)}{H^1_f(\mathbb{Q}_l, A^j_p)} \right),
\]
where $\mathrm{Loc}_l$ denotes the restriction of cohomology classes to the decomposition group at $l$ and the products run over all prime numbers.
\end{defn}
The $p$-part of Bloch-Kato's Tate-Shafarevich group $\Sha^{BK}_{\mathbb{Q}}(V^j_p)$ is defined in \cite[Definition 5.1]{BK} as follows.
\begin{defn}\label{Sha}
We define the $p$-part of Bloch-Kato's Tate-Shafarevich group for $V_p^j(=\Sym^j (V_p E))$ as 
\[
\Sha^{BK}_{\mathbb{Q}}(V^j_p):=\frac{H^1_f(\mathbb{Q}, A^j_p)}{\pi (H^1_f(\mathbb{Q}, V^j_p))},
\]
where the cohomology groups $H^1_f(\mathbb{Q}, A^j_p), H^1_f(\mathbb{Q}, V^j_p)$ are defined as in Definition \ref{Sel} and $\pi: H^1 (\mathbb{Q}, V^j_p) \rightarrow H^1 (\mathbb{Q}, A^j_p)$ is the canonical homomorphism induced by the natural map $\pi: V^j_p \rightarrow A^j_p$. In other words, $\Sha^{BK}_{\mathbb{Q}}(V^j_p)$ is defined by the exact sequence
\[
0 \rightarrow \pi (H^1_f(\mathbb{Q}, V^j_p)) \rightarrow H^1_f(\mathbb{Q}, A^j_p) \rightarrow \Sha^{BK}_{\mathbb{Q}}(V^j_p) \rightarrow 0.
\]
\end{defn}

Now we give a sketch of the proof of Theorem \ref{main}. In the following argument, we assume that the conditions $(a^{\prime}), (b^{\prime}), (c^{\prime})$ and $(d^{\prime})$ in Theorem \ref{main} hold.\\

\noindent\textbf{(Step1)}
We show the restriction map 
\[
\mathrm{Res}_{K/\mathbb{Q}}: H^1(\mathbb{Q}, \Sym^j E[p]) \rightarrow H^1(K, \Sym^j E[p])^G
\]
is an isomorphism where $G$ denotes $\Gal(K/\mathbb{Q})$. 
\vspace{2mm}

For a number field $F$, we define the unramified cohomology group $H_{\mathrm{ur}}^1(F, \Sym^j E[p])$ as the subgroup of cohomology classes in $H^1(F, \Sym^j E[p])$ which are trivial on the inertia group at every place of $F$. Assuming the claim in (Step1), the restriction $\mathrm{Res}_{K/\mathbb{Q}}$ induces an injective homomorphism between unramified cohomology groups
\begin{eqnarray*}
\mathrm{Res}_{K/\mathbb{Q}}:  H_{\mathrm{ur}}^1(\mathbb{Q}, \Sym^j E[p]) \rightarrow H_{\mathrm{ur}}^1(K, \Sym^j E[p])^G.
\end{eqnarray*}
Using class field theory, we have $H_{\mathrm{ur}}^1(K, \Sym^j E[p])^G=\mathrm{Hom}_G(\Cl_K \otimes \mathbb{F}_p, \Sym^j E[p])$. Every nontrivial homomorphism in $\mathrm{Hom}_G(\Cl_K\otimes \mathbb{F}_p, \Sym^j E[p])$ is surjective since $\Sym^j E[p]$ is irreducible. Thus the condition $H_{\mathrm{ur}}^1(\mathbb{Q}, \Sym^j E[p])\neq 0$ implies that $\Cl_K\otimes \mathbb{F}_p$ has $\Sym^j E[p]$ as its quotient $G$-module. We will construct nontrivial elements in $H_{\mathrm{ur}}^1(\mathbb{Q}, \Sym^j E[p])$ using Bloch-Kato's Selmer group in the succeeding steps.
\\

\noindent\textbf{(Step2)}
We show that the image of $H_f^1(\mathbb{Q}, \Sym^j E[p])$ in $H^1(\mathbb{Q}_l^{\mathrm{ur}}, \Sym^j E[p])$ is zero for any prime number $l\neq p$.
\vspace{2mm}

Here the cohomology group $H_f^1(\mathbb{Q}, \Sym^j E[p])$ is defined as follows. 
We have an exact sequence
\[
0 \rightarrow \Sym^j E[p] \xrightarrow{\iota} A^j_p \xrightarrow{\times p} A^j_p \rightarrow 0
\]
 from which we obtain an exact sequence
\[
0 \rightarrow \frac{H^0(\mathbb{Q}, A^j_p)}{pH^0(\mathbb{Q}, A^j_p)}\rightarrow H^1(\mathbb{Q}, \Sym^j E[p]) \xrightarrow{\iota} H^1(\mathbb{Q}, A^j_p)[p] \rightarrow 0, 
\]
where the map $\iota$ in the first exact sequence denotes the inclusion. We define the cohomology group $H_f^1(\mathbb{Q}, \Sym^j E[p])$ as the inverse image of the $p$-torsion part of Bloch-Kato's Selmer group $H_f^1(\mathbb{Q}, A^j_p)[p]$ under $\iota$. 
Assuming the claim in (Step2), for the restriction map 
\[
\mathrm{Res}^{\mathrm{ur}}_p: H_f^1(\mathbb{Q}, \Sym^j E[p]) \rightarrow H^1(\mathbb{Q}^{\mathrm{ur}}_p, \Sym^j E[p]),
\]
we have $\mathrm{Ker}(\mathrm{Res}^{\mathrm{ur}}_p) \subset H_{\mathrm{ur}}^1(\mathbb{Q}, \Sym^j E[p])$. Thus it suffices to show $\mathrm{Ker}(\mathrm{Res}^{\mathrm{ur}}_p) \neq 0$ to get the main theorem.\\

\noindent\textbf{(Step3)}
We study the image of $\mathrm{Res}^{\mathrm{ur}}_p$ and prove that $\dim_Fp (\mathrm{Im}(\mathrm{Res}^{\mathrm{ur}}_p)) \leqslant j$.
\vspace{2mm}

In Definition \ref{Sha}, we have an exact sequence
\[
0 \rightarrow \pi (H^1_f(\mathbb{Q}, V^j_p)) \rightarrow H^1_f(\mathbb{Q}, A^j_p) \rightarrow \Sha^{BK}_{\mathbb{Q}}(V^j_p) \rightarrow 0.
\]
Since the group $\pi (H^1_f(\mathbb{Q}, V^j_p))$ is $p$-divisible, the above homomorphism $ H^1_f(\mathbb{Q}, A^j_p)\rightarrow \Sha^{BK}_{\mathbb{Q}}(V^j_p)$ is still surjective when restricted on the $p$-torsion parts. Since $H_f^1(\mathbb{Q}, \Sym^j E[p])$ is defined as the inverse image of $H^1_f(\mathbb{Q}, A^j_p)[p]$ under the surjection $\iota$, we have a surjective map $H^1_f(\mathbb{Q}, \Sym^j E[p]) \twoheadrightarrow \Sha^{BK}_{\mathbb{Q}}(V^j_p)[p].$ So if we assume the condition $\dim_Fp(\Sha^{BK}_{\mathbb{Q}}(V^j_p)[p])\geqslant j+1$ in Theorem \ref{main}, then we have $\dim_Fp (H^1_f(\mathbb{Q}, \Sym^j E[p]))\geqslant j+1$. From the claim in (Step3), we have $\mathrm{Ker}(\mathrm{Res}^{\mathrm{ur}}_p) \neq 0$ and the theorem follows.
\begin{rem}\label{rem}
We can apply the above argument to more general $p$-adic Galois representations. For example, the representations attached to modular forms can be treated. Let $f$ be a normalized new eigen cusp form whose coefficients are in $\mathbb{Q}$ and level prime to $p$. For this modular form $f$, we have an associated integral $p$-adic Galois representation $\rho_{f, p}:G_{\mathbb{Q}} \rightarrow \mathrm{Aut}_{\mathbb{Z}_p}(T_{f, p}) \cong \mathrm{GL}_2(\mathbb{Z}_p)$, where $T_{f, p}$ is its representation space which is a free $\mathbb{Z}_p$ module of rank 2. Let $\overline{\rho}_{f, p}:G_{\mathbb{Q}} \rightarrow \mathrm{Aut}_{\mathbb{F}_p}(\overline{T}_{f, p}) \cong \mathrm{GL}_2(\mathbb{F}_p)$ be the mod $p$ reduction of $\rho_{f, p}$. We consider twists of these representations. For a square-free integer $D$ and $j \in \mathbb{Z}$, let $\rho_{f, p}(j, D):= \rho_{f, p} \otimes \chi_{D} \otimes \chi^j_{\mathrm{cyc}}:G_{\mathbb{Q}} \rightarrow \mathrm{Aut}_{\mathbb{Z}_p}(T_{f, p}(j, D)) \cong \mathrm{GL}_2(\mathbb{Z}_p)$ be a twist of $\rho_{f, p}$ by a quadratic Dirichlet character $\chi_{D}$ associated to $D$ and $\chi_{\mathrm{cyc}}^j$, where $\chi_{\mathrm{cyc}}$ is the $p$-adic cyclotomic character. Let $\overline{\rho}_{f, p}(j, D):G_{\mathbb{Q}} \rightarrow \mathrm{GL}_2(\mathbb{F}_p)$ be the mod $p$ reduction of $\rho_{f, p}(j, D)$ corresponding to $\overline{T}_{f, p}(j, D):=T_{f, p}(j, D) \otimes \mathbb{Z}/p$. For the $\mathbb{F}_p$-valued representation $\overline{\rho}_{f, p}(j, D)$, we have a number field $K_{f, p, (j, D)}$ corresponding to the kernel of $\overline{\rho}_{f, p}(j, D)$. Let $\mathrm{Cl}_{f, p, (j, D)}$ be the ideal class group of $K_{f, p, (j, D)}$. We can apply the same argument in this article to the $p$-adic representation $V_{f, p}(j, D):=T_{f, p}(j, D)\otimes \mathbb{Q}_p$ under similar assumptions to those in Theorem \ref{main}. Then we can deduce that the condition $\dim_Fp(\Sha^{BK}_{\mathbb{Q}}(V_{f, p}(j, D))[p])\geqslant 2$ implies that the $\mathbb{F}_p$-representation $\mathrm{Cl}_{f, p, (j, D)}\otimes \mathbb{F}_p$ of $\Gal(K_{f ,p, (j, D)}/\mathbb{Q})$ has $\overline{T}_{f, p}(j, D)$ as its quotient representation. Assuming the Bloch-Kato conjecture, we can make some numerical examples of the above result with some calculations of special values of $L$-functions attached to $f$ and various $\chi_{D}$. We will describe the details of this result in our forthcoming paper.
\end{rem}
\section{Injectivity of the restriction map}
In this section, we prove the claim in (Step1) in the previous section.
\begin{prop}\label{inj}
Suppose the representation $\rho_{E, p}$ is surjective. Then the restriction map 
\[
\mathrm{Res}_{K/\mathbb{Q}}: H^1(\mathbb{Q}, \Sym^j E[p]) \rightarrow H^1(K, \Sym^j E[p])^G
\]
is an isomorphism.
\end{prop}
(Proof of Proposition \ref{inj}) 

It suffices to show that $H^1(G, \Sym^j E[p])=H^2(G, \Sym^j E[p])=0$. We use the following lemma.
\begin{lem}\label{vanishmentofcoh}
Let $G$ be a finite group and $V$ a finite dimensional representation of $G$ over a field $F$ of characteristic $p$. If there is a normal subgroup $H$ of $G$ such that 
\begin{itemize}
\item[(1)] $\#H$ is prime to  $p$
\item[(2)] $V^{H}=0$
\end{itemize}
then  $H^i(G, V)=0$ for all $i\geqslant 0$.
\end{lem}
(Proof of lemma \ref{vanishmentofcoh})

The condition (2) implies $H^0(G, V)=0$. We have the inflation-restriction exact sequence
\[
0 \rightarrow H^1(G/H, V^H) \xrightarrow{\textrm{inf}} H^1(G, V) \xrightarrow{\textrm{res}} H^1(H, V)^G.
\]
Since $V^H=0$ and $\#H$ is prime to $p$, the first and the third term in the above sequence are 0 and we get $H^1(G, V)=0$. Since $H^1(H, V)=0$, we also have the inflation-restriction exact sequences for the cohomology groups of higher degree inductively to get $H^i(G, V)=0$ for $i\geqslant 0$. \hfill$\square$

\vspace{2mm}
 We go back to the proof of the proposition. Since we assume $1\leqslant j\leqslant p-2$,  there is an element $c \in \mathbb{F}_p^{\times}$ such that $c^j \neq 1$. The central element $cI$ acts on $E[p]$ by multiplication by $c$, here $I$ denotes the unit matrix in $\textrm{GL}_2(\mathbb{F}_p)$. Then $cI$ acts on $\Sym^j E[p]$ by multiplication by $c^j (\neq 1)$. Let $C$ be the subgroup of $G$ generated by $cI$. Since $cI$ is a central element and $c \in \mathbb{F}^{\times}_p$, $C$ is a normal subgroup of $G$ and $\# C$ is prime to $p$. So the subgroup $C$ satisfies the conditions $(1), (2)$ of Lemma \ref{vanishmentofcoh}, then we have $H^i(G, \Sym^j E[p])=0$ for $i\geqslant 0$ and $1\leqslant j \leqslant p-2$. Hence the injectivity of $\mathrm{Res}_{K/\mathbb{Q}}: H^1(\mathbb{Q}, \Sym^j E[p]) \rightarrow H^1(K, \Sym^j E[p])^G
$ follows. \hfill$\square$
 
\section{The cohomology group $H_f^1(\mathbb{Q}, \Sym^j E[p])$}
Next we show the claim in (Step2) in Section 2.
\begin{prop}\label{unram}
For a prime $l\neq p$, suppose $v_{l}(j(E))$ is prime to $p$ when $E$ has potentially multiplicative reduction at $l$. Then the elements in $H_f^1(\mathbb{Q}, \Sym^j E[p])$ are unramified outside $p$.
\end{prop}
(Proof of Proposition \ref{unram})

 Since $p$ is an odd prime, any elements in $H^1(\mathbb{Q}, \Sym^j E[p])$ are unramified at the infinite place of $\mathbb{Q}$ automatically.
 
 For every prime number $l \neq p$, we have the following commutative diagram
\[
  \xymatrix{
  &    & H_f^1(\mathbb{Q},\Sym^j E[p]) \ar[r]^{\iota}\ar[d]^{\mathrm{Res}_{l}^{\mathrm{ur}}} & H_f^1(\mathbb{Q}, A^j_p)[p]\ar[r]^{} \ar[d]^{\mathrm{res}_{l}^{\mathrm{ur}}}& 0\\
 0 \ar[r]^{}&   \displaystyle  \frac{H^0(\mathbb{Q}^{\mathrm{ur}}_l, A^j_p)}{pH^0(\mathbb{Q}^{\mathrm{ur}}_l, A^j_p)} \ar[r] &  H^1(\mathbb{Q}^{\mathrm{ur}}_l, \Sym^j E[p]) \ar[r]^{\iota} & H^1(\mathbb{Q}^{\mathrm{ur}}_l, A^j_p)[p] \ar[r]^{} & 0.
  }
\]
Here $\mathrm{Res}_{l}^{\mathrm{ur}}$ denotes the restriction of cohomology classes to the inertia at $l$ and $\iota$ is the homomorphism induced by the inclusion $\Sym^j E[p] \hookrightarrow A^j_p$. What we have to show is that for any cohomology classes $c  \in H_f^1(\mathbb{Q}, \Sym^j E[p])$, we have $\mathrm{Res}_{l}^{\mathrm{ur}}(c)=0$. So it suffices to show $\frac{H^0(\mathbb{Q}_{l}^{\mathrm{ur}}, A^j_p)}{pH^0(\mathbb{Q}_{l}^{\mathrm{ur}}, A^j_p)}=0$.\vspace{1mm}

\noindent(Case 1) $E$ has good reduction at $l$.

\vspace{1mm}
In this case, $\Sym^j E[p^n]$ is unramified at $l$ for any positive integer $n$. So we have $H^0(\mathbb{Q}_{l}^{\mathrm{ur}}, A^j_p)=(\mathbb{Q}_p/\mathbb{Z}_p)^{\oplus(j+1)}$ to get $\frac{H^0(\mathbb{Q}_{l}^{\mathrm{ur}}, A^j_p)}{pH^0(\mathbb{Q}_{l}^{\mathrm{ur}}, A^j_p)}=0$.\vspace{1mm}

\noindent(Case 2) $E$ has split multiplicative reduction at $l$.

\vspace{1mm}
In this case, using the result of Tate, we have an isomorphism $E(\overline{\mathbb{Q}_{l}}) \cong \overline{\mathbb{Q}_{l}}^{\times}/\langle q \rangle$ as $G_{\mathbb{Q}_{l}}$-modules, here $q$ is the Tate period for $E$ in $\mathbb{Q}_{l}$. Then for a positive integer $n$, the group of $p^n$-torsion points $E[p^n]$ is isomorphic to a free $\mathbb{Z}/p^n$-module generated by $\zeta_{p^n}$ and $\sqrt[p^n]{q}$, where $\zeta_{p^n}$ and $\sqrt[p^n]{q}$ denote a primitive $p^n$-th root of unity and a $p^n$-th root of $q$ respectively. So with respect to this basis, $G_{\mathbb{Q}_{l}^{\mathrm{ur}}}$ acts on $E[p^n]$ via
\[
 \left(
 \begin{array}{cc}
 1&\tau_{q, n}\\
 0&1
 \end{array}
 \right),
\]
where $\tau_{q, n}:G_{\mathbb{Q}^{\mathrm{ur}}_{l}} \rightarrow \mathbb{Z}/{p^n}$ is the map defined by $
g(\sqrt[p^n]{q})=\sqrt[p^n]{q}\cdot \zeta_{p^n}^{\tau_{q, n}(g)}\ \ (g \in G_{\mathbb{Q}^{\mathrm{ur}}_{l}}).
$
For a positive integer $n$, we compute $H^0(\mathbb{Q}_{l}^{\mathrm{ur}}, \Sym^j E[p^n])$ explicitly. First, we fix a basis of $\Sym^j E[p^n]$ over $\mathbb{Z}/{p^n}$ as 
\[
u_0:=\zeta_{p^n}^{\otimes j},\ u_1:=\zeta_{p^n}^{\otimes j-1} \otimes \sqrt[p^n]{q},\ \ldots, u_i:=\zeta_{p^n}^{\otimes j-i} \otimes \sqrt[p^n]{q}^{\otimes i},\ \ldots, \ u_j:=\sqrt[p^n]{q}^{\otimes j}.
\]
\begin{lem}\label{comp}
For an element $x:=a_0 u_0+a_1 u_1+\cdots +a_j u_j \in \Sym^j E[p^n]$ with $a_i \in \mathbb{Z}/{p^n}$, the condition $x \in \Sym^j E[p^n]^{G_{\mathbb{Q}_{l}^{\mathrm{ur}}}}$ is equivalent to the condition
\[
a_0 \in \mathbb{Z}/{p^n},\ a_1\tau_{q, n}(g)=a_2\tau_{q, n}(g)=\cdots=a_j\tau_{q, n}(g)=0\ \ \text{in}\ \mathbb{Z}/{p^n}\ \ (\forall g \in G_{\mathbb{Q}_{l}^{\mathrm{ur}}}).
\]
\end{lem}
(Proof of Lemma \ref{comp})\\
This can be proved by an explicit calculation. For $i$ with $0\leqslant i\leqslant j$ and $g \in G_{\mathbb{Q}_{l}^{\mathrm{ur}}}$, we have 
\begin{eqnarray*}
g(u_i)=g(\zeta_{p^n}^{\otimes j-i} \otimes \sqrt[p^n]{q}^{\otimes i})&=& \zeta_{p^n}^{\otimes j-i}\otimes (\tau_{q, n}(g)\zeta_{p^n}+\sqrt[p^n]{q})^{\otimes i}=\sum_{k=0}^{i}\binom{i}{k}\tau_{q, n}(g)^k u_{i-k}.
\end{eqnarray*}
So for $x:=a_0 u_0+a_1 u_1+\cdots +a_j u_j \in \Sym^j E[p^n]$, we can deduce that the condition $x \in \Sym^j E[p^n]^{G_{\mathbb{Q}_{l}^{\mathrm{ur}}}}$ is equivalent to the condition
\begin{eqnarray}\label{induction}
\sum_{k=i+1}^{j}\binom{k}{i}a_k \tau_{q, n}(g)^{k-i}=0\ \ in\ \mathbb{Z}/{p^n}\ \ (0\leqslant i\leqslant j-1)
\end{eqnarray}
for all $g \in G_{\mathbb{Q}_{l}^{\mathrm{ur}}}$. When $i=j-1$, we have $\binom{j}{j-1}a_j\tau_{q, n}(g)=0$ from the equation (\ref{induction}) to get $a_j\tau_{q, n}(g)=0$ since $j\leqslant p-2$. When $i=j-2$, we have 
\[
\binom{j-1}{j-2}a_{j-1}\tau_{q, n}(g)+\binom{j}{j-2}a_j \tau_{q, n}(g)^2=0
\]
from (\ref{induction}). Since $a_j\tau_{q, n}(g)=0$ and $j\leqslant p-2$, we also have $a_{j-1}\tau_{q, n}(g)=0$. By backward induction on $i$, we can get $a_i\tau_{q, n}(g)=0$ for $1\leqslant i \leqslant j$. \hfill$\square$
\\

We go back to the proof of Proposition \ref{unram}. If there is an element $g \in G_{\mathbb{Q}_{l}^{\mathrm{ur}}}$ such that $\tau_{q, n}(g) \in (\mathbb{Z}/{p^n})^{\times}$, then we have $a_1=a_2=\cdots=a_j=0$ and $\Sym^j E[p^n]^{G_{\mathbb{Q}_{l}^{\mathrm{ur}}}}\cong \mathbb{Z}/{p^n}$ for any $n$ from Lemma \ref{comp}. Thus we have $H^0(\mathbb{Q}_{l}^{\mathrm{ur}}, A^j_p)\cong \mathbb{Q}_p/\mathbb{Z}_p$ and $\frac{H^0(\mathbb{Q}_{l}^{\mathrm{ur}}, A^j_p)}{pH^0(\mathbb{Q}_{l}^{\mathrm{ur}}, A^j_p)}=0$. So we will show in the following that there exists such $g \in G_{\mathbb{Q}_{l}^{\mathrm{ur}}}$ under the assumptions in Proposition \ref{unram}. It is a well-known fact that if an elliptic curve $E$ over $\mathbb{Q}_{l}$ has split multiplicative reduction, then $v_{l}(j(E))=-v_{l}(q)$ where $v_{l}$ is the normalized $l$-adic valuation on $\mathbb{Q}_{l}$ and $q$ is the Tate period for $E$. Since we assume $v_{l}(j(E))$ is prime to $p$, $q$ is not a $p$-th power in $\mathbb{Q}_{l}$. So we have $\sqrt[p^n]{q} \notin \mathbb{Q}_{l}$ for any $n \in \mathbb{Z}_{>0}$ and $\sqrt[p^n]{q} \notin \mathbb{Q}_{l}^{\mathrm{ur}}$ since $v_{l}(q)>0$. Then there exists $g\in G_{\mathbb{Q}_{l}^{\mathrm{ur}}}$ such that $\tau_{q, n}(g) \in (\mathbb{Z}/{p^n})^{\times}$ for any $n$ and the proposition follows in this case.\vspace{1mm}

\noindent(Case 3) $E$ has non-split multiplicative reduction at $l$.

\vspace{1mm}
In this case, $E$ has split multiplicative reduction at $l$ over the unramified quadratic extension $F$ of $\mathbb{Q}_{l}$. So we can imitate  the argument in the (Case 2) over $F$ to get the desired result.\vspace{1mm}

\noindent(Case 4) $E$ has additive potentially multiplicative reduction at $l$.

\vspace{1mm}
In this case, $E$ has split multiplicative reduction at a prime above $l$ over a ramified quadratic extension $L$ of $\mathbb{Q}_{l}$. So there exists some quadratic twist $E^{\prime}$ of the elliptic curve $E$, which has split multiplicative reduction at $l$ and for each positive integer $n$, as $G_{\mathbb{Q}_{l}}$-modules, we have
\[
E[p^n] \cong E^{\prime}[p^n] \otimes \chi
\]
here $\chi$ denotes the ramified quadratic character corresponds to $L$. Taking the Tate period $q$ for $E^{\prime}$, for some suitable basis $\{v_1, v_2\}$, we know that the action of $G_{\mathbb{Q}_{l}^{\mathrm{ur}}}$ on $E[p^n]$ is of the form:
\[
 \left(
 \begin{array}{cc}
 1&\tau_{q, n}\\
 0&1
 \end{array}
 \right)\otimes \chi
\]
as in the argument in (Case 2). We again fix a basis of $\Sym^j E[p^n]$ over $\mathbb{Z}/p^n\mathbb{Z}$ as 
\[
u_0:=v_1^{\otimes j},\ u_1:=v_1^{\otimes j-1} \otimes v_2,\ \ldots, u_i:=v_1^{\otimes j-i} \otimes v_2^{\otimes i},\ \ldots, \ u_j:=v_2^{\otimes j}.
\]
 For an element $x:=a_0u_0+a_1u_1+\cdots+a_j u_j \in \Sym^j E[p^n]$ with $a_i \in \mathbb{Z}/p^n\mathbb{Z}$, we can show that the condition $x \in \Sym^j E[p^n]^{G_{L\cdot\mathbb{Q}_{l}^{\mathrm{ur}}}}$ is equivalent to the condition
\[
a_0 \in \mathbb{Z}/p^n\mathbb{Z},\ a_1\tau_{q, n}(g)=a_2\tau_{q, n}(g)=\cdots=a_j\tau_{q, n}(g)=0\ \ \textrm{in}\ \mathbb{Z}/{p^n}\ \ (\forall g \in G_{L\cdot\mathbb{Q}_{l}^{\mathrm{ur}}})
\]
by the same calculation as in the proof of Lemma \ref{comp} since $\chi$ is trivial on $G_{L\cdot\mathbb{Q}_{l}^{\mathrm{ur}}}$.
 On the other hand, we have $\sqrt[p^n]{q} \notin \mathbb{Q}_{l}^{\mathrm{ur}}$ for any positive integer $n$ by the assumption that $p$ does not divide $v_{l}(q)=-v_{l}(j(E^{\prime}))=-v_{l}(j(E))$. We know $\sqrt[p^n]{q}$ is also not contained in $L\cdot \mathbb{Q}_{l}^{\mathrm{ur}}$ because $L\cdot \mathbb{Q}_{l}^{\mathrm{ur}}/\mathbb{Q}_{l}^{\mathrm{ur}}$ is a quadratic extension and $p \neq 2$. So there exists $g \in G_{L\cdot\mathbb{Q}_{l}^{\mathrm{ur}}}$ such that $\tau_{q, n}(g) \in (\mathbb{Z}/{p^n})^{\times}$ for every $n$ and we get $a_1=a_2=\ldots=a_j=0$ as in the argument in (Case 2). Thus we get $(\Sym^j E[p^n])^{G_{L\cdot\mathbb{Q}_{l}^{\mathrm{ur}}}}=\mathbb{Z}/{p^n}\cdot u_0$ and $(\Sym^j E[p^n])^{G_{\mathbb{Q}_{l}^{\mathrm{ur}}}}=(\mathbb{Z}/{p^n}\cdot u_0)^{{\rm{Gal}}(L/\mathbb{Q}_{l})}$. Let $\tau$ denote a generator of the Galois group ${\rm{Gal}}(L/\mathbb{Q}_{l})$. Then $\tau(u_0)=\chi(\tau)^j u_0=(-1)^j u_0$. So we have
\[
\Sym^j E[p^n]^{G_{\mathbb{Q}^{\mathrm{ur}}_{l}}} = \left\{\begin{array}{ll}
0 & (\text{$j$ is odd}) \\  
\mathbb{Z}/{p^n}\cdot u_0 & (\text{$j$ is even}).
\end{array}\right.
\]
We have $\frac{H^0(\mathbb{Q}_{l}^{\mathrm{ur}}, A^j_p)}{pH^0(\mathbb{Q}_{l}^{\mathrm{ur}}, A^j_p)}=0$ in both cases. \vspace{1mm} 

\noindent(Case 5) $E$ has additive potentially good reduction at $l$.

Let $\rho:G_{\mathbb{Q}} \rightarrow \mathrm{GL}_2(\mathbb{Z}_p)$ be the representation associated to the integral $p$-adic Tate module of $E$. It is a well-known fact that when $E$ has potentially good reduction at $l$, then $\#\rho(G_{\mathbb{Q}_{l}^{\mathrm{ur}}})$ is finite and its possible prime divisors are only 2 and 3. See for example, \cite[Section 3.3]{FK}. So we get $p\nmid \#\rho(G_{\mathbb{Q}_{l}^{\mathrm{ur}}})$ since we assume $p\geqslant5$. Let $\rho_j : G_{\mathbb{Q}} \rightarrow \mathrm{GL}_{j+1}(\mathbb{Z}_p)$ be the representation attached to $T_p^j$. Then we have $\mathrm{Ker}(\rho) \subset \mathrm{Ker}(\rho_j)$ to get a natural surjection $\rho(G_{\mathbb{Q}_{l}^{\mathrm{ur}}}) \twoheadrightarrow \rho_j(G_{\mathbb{Q}_{l}^{\mathrm{ur}}})$. Thus we also get $p \nmid \#\rho_j(G_{\mathbb{Q}_{l}^{\mathrm{ur}}})$. This implies that there is an open normal subgroup $U$ of $G_{\mathbb{Q}_{l}^{\mathrm{ur}}}$ such that $\rho_j(U)=0$ and $[G_{\mathbb{Q}_{l}^{\mathrm{ur}}} : U]$ is prime to $p$. Then we have the inflation-restriction exact sequence
\[
0 \rightarrow H^1(G_{\mathbb{Q}_{l}^{\mathrm{ur}}}/U, (T^j_p)^U) \rightarrow H^1(G_{\mathbb{Q}_{l}^{\mathrm{ur}}}, T^j_p) \rightarrow H^1(U, T^j_p)^{G_{\mathbb{Q}_{l}^{\mathrm{ur}}}/U} \rightarrow H^2(G_{\mathbb{Q}_{l}^{\mathrm{ur}}}/U, (T^j_p)^U).
\]
Since the order of $G_{\mathbb{Q}_{l}^{\mathrm{ur}}}/U$ is prime to $p$, we have $H^i(G_{\mathbb{Q}_{l}^{\mathrm{ur}}}/U, (T^j_p)^U)=0$ for $i=1, 2$ and obtain an isomorphism
\[
H^1(G_{\mathbb{Q}_{l}^{\mathrm{ur}}}, T^j_p) \cong H^1(U, T^j_p)^{G_{\mathbb{Q}_{l}^{\mathrm{ur}}}/U}
\]
induced by the restriction map. We know that $U$ acts trivially on $T^j_p$ to get $H^1(U, T^j_p)=\mathrm{Hom}(U, T^j_p)$. Since $T^j_p$ is torsion-free, this group $\mathrm{Hom}(U, T^j_p)$ and of course its subgroup $H^1(U, T^j_p)^{G_{\mathbb{Q}_{l}^{\mathrm{ur}}}/U}$ are torsion-free. On the other hand, we have an exact sequence $0 \rightarrow T_p^j \rightarrow V_p^j \rightarrow A_p^j \rightarrow 0$ from which we also have an exact sequence
\[
0 \rightarrow (T^j_p)^{G_{\mathbb{Q}_{l}^{\mathrm{ur}}}}\otimes \mathbb{Q}_p/\mathbb{Z}_p \rightarrow (A^j_p)^{G_{\mathbb{Q}_{l}^{\mathrm{ur}}}} \rightarrow H^1(G_{\mathbb{Q}_{l}^{\mathrm{ur}}}, T^j_p)[p^{\infty}] \rightarrow 0.
\]
Since $H^1(G_{\mathbb{Q}_{l}^{\mathrm{ur}}}, T^j_p)[p^{\infty}]=0$ from the above argument, we have $(T^j_p)^{G_{\mathbb{Q}_{l}^{\mathrm{ur}}}}\otimes \mathbb{Q}_p/\mathbb{Z}_p \cong (A^j_p)^{G_{\mathbb{Q}_{l}^{\mathrm{ur}}}}$ and $(A^j_p)^{G_{\mathbb{Q}_{l}^{\mathrm{ur}}}}=H^0(G_{\mathbb{Q}_{l}^{\mathrm{ur}}}, A^j_p)$ is divisible. Thus $\frac{H^0(\mathbb{Q}_{l}^{\mathrm{ur}}, A^j_p)}{pH^0(\mathbb{Q}_{l}^{\mathrm{ur}}, A^j_p)}=0$. We have proved the Proposition \ref{unram} in all cases. \hfill$\square$
\section{The image of the restriction map $\mathrm{Res}_p^{\mathrm{ur}}$}
We finally prove the following proposition which is the claim in (Step 3).
\begin{prop}\label{dimofResp}
For the restriction map 
\[
\mathrm{Res}_p^{\mathrm{ur}}:H_f^1(\mathbb{Q}, \Sym^j E[p]) \rightarrow H^1(\mathbb{Q}^{\mathrm{ur}}_p, \Sym^j E[p]),
\]
we have
\[
\dim_Fp (\mathrm{Im}(\mathrm{Res}^{\mathrm{ur}}_p)) \leqslant j.
\]
\end{prop}
For the above restriction map, we have a decomposition
\[
\mathrm{Res}^{\mathrm{ur}}_p:H_f^1(\mathbb{Q},\Sym^j E[p])\xrightarrow{{\mathrm{Loc}_p}} H^1(\mathbb{Q}_{p}, \Sym^j E[p]) \xrightarrow{\mathrm{Res}_{\mathbb{Q}_p^{\mathrm{ur}}/\mathbb{Q}_p}} H^1(\mathbb{Q}^{\mathrm{ur}}_{p}, \Sym^j E[p]).
\]
Here the homomorphism $\mathrm{Loc}_p$ is the restriction of cohomology classes to the decomposition group at $p$, and $\mathrm{Res}_{\mathbb{Q}_p^{\mathrm{ur}}/\mathbb{Q}_p}$ is the restriction of them to the inertia group at $p$. So first we study the image of $\mathrm{Loc}_p$. We have the following commutative diagram
\[
  \xymatrix{
  &    & H_f^1(\mathbb{Q},\Sym^j E[p]) \ar[r]^{\iota}\ar[d]^{\mathrm{Loc}_p} & H_f^1(\mathbb{Q}, A^j_p)[p]\ar[r]^{} \ar[d]^{\mathrm{loc}_{p}}& 0\\
 0 \ar[r]^{}&   \displaystyle  \frac{H^0(\mathbb{Q}_p, A^j_p)}{pH^0(\mathbb{Q}_p, A^j_p)} \ar[r] &  H^1(\mathbb{Q}_{p}, \Sym^j E[p]) \ar[r]^{\iota} & H^1(\mathbb{Q}_{p}, A^j_p)[p] \ar[r]^{} & 0. 
  }
\]
So we have $\mathrm{Im}(\mathrm{Loc}_p) \subset \iota^{-1}(H_f^1(\mathbb{Q}_{p}, A^j_p)[p])$ and there is an exact sequence
\begin{eqnarray}\label{locexact}
0 \rightarrow  \displaystyle  \frac{H^0(\mathbb{Q}_p, A^j_p)}{pH^0(\mathbb{Q}_p, A^j_p)} \rightarrow \iota^{-1}(H_f^1(\mathbb{Q}_{p}, A^j_p)[p]) \xrightarrow{\iota}  H_f^1(\mathbb{Q}_{p}, A^j_p)[p] \rightarrow 0.
\end{eqnarray}
 So we have an inequality
 \begin{eqnarray}\label{fundineq}
 \dim_Fp(\mathrm{Im}(\mathrm{Loc}_p)) \leqslant \displaystyle \dim_Fp\left(\frac{H^0(\mathbb{Q}_p, A^j_p)}{pH^0(\mathbb{Q}_p, A^j_p)}\right)+\dim_Fp\left( H_f^1(\mathbb{Q}_{p}, A^j_p)[p] \right).
 \end{eqnarray}
 The dimension of $ H_f^1(\mathbb{Q}_{p}, A^j_p)[p]$ can be computed by $p$-adic Hodge theory. We use the following fact in \cite[Section 9.2.2]{TQ}.
 \begin{prop}
 Let $V$ be a $p$-adic representation of $G_{\mathbb{Q}_p}$ and \\ $\mathbf{D}_{\mathrm{dR}}(V):=(V\otimes \mathbf{B}_{\mathrm{dR}})^{G_{\mathbb{Q}_p}}, 
 \mathbf{D}^{+}_{\mathrm{dR}}(V):=(V\otimes \mathbf{B}^{+}_{\mathrm{dR}})^{G_{\mathbb{Q}_p}}$, where $\mathbf{B}_{\mathrm{dR}}$ is the Fontaine's de Rham period ring. If $V$ is a de Rham representation, then
 \begin{eqnarray}\label{hodge}
 \mathrm{dim}_{\mathbb{Q}_p}(H_f^1(\mathbb{Q}_{p}, V))=  \mathrm{dim}_{\mathbb{Q}_p}(\mathbf{D}_{\mathrm{dR}}(V)/ \mathbf{D}^{+}_{\mathrm{dR}}(V))+\mathrm{dim}_{\mathbb{Q}_p}H^0(\mathbb{Q}_p, V).
 \end{eqnarray}
 \end{prop}
The $p$-adic representation $V_p E$ is crystalline with Hodge-Tate weight \{0,1\} because of the assumption that $E$ has good reduction at $p$. Since the functor $\mathbf{D}_{\mathrm{dR}}$ is compatible with taking symmetric powers, we can compute that $\mathrm{dim}_{\mathbb{Q}_p}(\mathbf{D}_{\mathrm{dR}}(V_p^j)/ \mathbf{D}^{+}_{\mathrm{dR}}(V_p^j))=(j+1)-1=j$. By the equality (\ref{hodge}) and the definition of $H_f^1(\mathbb{Q}_{p}, A^j_p)$, we have
 \begin{eqnarray}\label{ineq2}
 \dim_Fp\left( H_f^1(\mathbb{Q}_{p}, A^j_p)[p]\right) = j + \mathrm{dim}_{\mathbb{Q}_p}H^0(\mathbb{Q}_p, V).
 \end{eqnarray}
 From (\ref{fundineq}), (\ref{ineq2}), we have
 \begin{eqnarray}\label{ineq}
  \dim_Fp(\mathrm{Im}(\mathrm{Loc}_p)) \leqslant \displaystyle \dim_Fp\left(\frac{H^0(\mathbb{Q}_p, A^j_p)}{pH^0(\mathbb{Q}_p, A^j_p)}\right)+j + \textrm{dim}_{\mathbb{Q}_p}H^0(\mathbb{Q}_p, V).
 \end{eqnarray}
In the following, we compute the first term and the third term of the right-hand side of (\ref{ineq}).
\begin{prop}\label{dimssup}
Suppose that $E$ has good supersingular reduction at $p$. Then
 \[
 \displaystyle \dim_Fp\left(\frac{H^0(\mathbb{Q}_p, A^j_p)}{pH^0(\mathbb{Q}_p, A^j_p)}\right)=\mathrm{dim}_{\mathbb{Q}_p}H^0(\mathbb{Q}_p, V)=0.
 \]
\end{prop}
(Proof of Proposition \ref{dimssup})

Since the computations of $\dim_Fp \frac{H^0(\mathbb{Q}_p, A^j_p)}{pH^0(\mathbb{Q}_p, A^j_p)}$ and $\mathrm{dim}_{\mathbb{Q}_p}H^0(\mathbb{Q}_p, V^j_p)$ are very similar, we only describe precise computations for $\dim_Fp \frac{H^0(\mathbb{Q}_p, A^j_p)}{pH^0(\mathbb{Q}_p, A^j_p)}$.
 
 If $E$ has good supersingular reduction at $p$, for every positive integer $n$, $E[p^n]$ is isomorphic to the group of the $p^n$-torsion points of the Lubin-Tate formal group associated to a prime $-p$ over an unramified quadratic extension $F$ of $\mathbb{Q}_p$ (\cite[Proposition 8.6]{K}). So $E[p^n]$ is a free $\mathcal{O}_{F}/p^n$-module of rank 1 and we take its basis $z_n$ over $\mathcal{O}_{F}/p^n$, where $\mathcal{O}_{F}$ denotes the ring of integers of $F$. The Galois group $G_{\mathbb{Q}^{\mathrm{ur}}_p}$ acts on $z_n$ via the character $\overline{\chi_{LT}}:=\chi_{LT} \mod p^n$. Here $\chi_{LT}:G_{\mathbb{Q}^{\mathrm{ur}}_p} \twoheadrightarrow \mathcal{O}_{F}^{\times}$ is the Lubin-Tate character associated to the prime element $-p$ of $F$. We take $y_n:=z_n^{\otimes j}$ as a basis of $\Sym^j E[p^n]$ over $\mathcal{O}_{F}/p^n$. For an element $x=ay_n\ (a\in \mathcal{O}_{F}/p^n)$ in $\Sym^j E[p^n]$, the condition that $x \in \Sym^j E[p^n]^{G_{\mathbb{Q}^{\mathrm{ur}}_p}}$ is equivalent to the condition that $a(\overline{\chi_{LT}}^{j}(g)-1)=0$ for all $g \in G_{\mathbb{Q}^{\mathrm{ur}}_p}$. Since we assume $j<p-1$, there exists $g\in G_{\mathbb{Q}^{\mathrm{ur}}_p}$ such that $\overline{\chi_{LT}}^{j}(g)-1 \in (\mathcal{O}_{F}/p^n)^{\times}$ and we get $a=0$. Hence we have $\Sym^j E[p^n]^{G_{\mathbb{Q}_p}}=\Sym^j E[p^n]^{G_{\mathbb{Q}^{\mathrm{ur}}_p}}=0$ for all $n$ and $H^0(\mathbb{Q}_p, A^j_p)=\Sym^j E[p^{\infty}]^{G_{\mathbb{Q}_p}}=0$.\\
 \hfill$\square$

Here we introduce some notations for the good ordinary reduction case. If $E$ has good ordinary reduction at $p$, we have an exact sequence 
\[
0 \rightarrow T_p \widehat{E} \rightarrow T_p E \rightarrow T_p \widetilde{E}_p \rightarrow 0,
\]
where $ \widetilde{E}_p$ is the mod $p$ reduction of the curve $E$ and $ \widehat{E}$ is the kernel of the reduction. We take a basis $\{v_1, v_2\}$ of  $T_p E$ as a $\mathbb{Z}_p$-module such that $T_p \widehat{E} = \mathbb{Z}_p v_1$ and we have the representation $\rho_{E} : G_{\mathbb{Q}_p} \rightarrow \mathrm{GL}_2(\mathbb{Z}_p)$ with respect to this basis. For each positive integer $n$, \{$v_1$ mod $p^n$, $v_2$ mod $p^n$\} form a basis of the free $\mathbb{Z}/p^n$-module $E[p^n]$ and this basis yields the representation $\rho_{E, p^n}: G_{\mathbb{Q}_p} \rightarrow \mathrm{GL}_2(\mathbb{Z}/p^n)$. The action of $g \in G_{\mathbb{Q}_p}$ on $T_p E$ can be written as the matrix
\begin{eqnarray}\label{action}
 \left(
 \begin{array}{cc}
 \chi_{\mathrm{cyc}}(g)\psi^{-1}(g) &u(g) \\
 0&\psi(g)
 \end{array}
 \right).
\end{eqnarray}
Here $\chi_{\mathrm{cyc}}$ denotes the $p$-adic cyclotomic character, $\psi$ is the unramified character determined by the action of $G_{\mathbb{Q}_p}$ on $T_p \widetilde{E}_p$, and $u(g) \in \mathbb{Z}_p$.
Also the action of $g \in G_{\mathbb{Q}_p}$ on $E[p^n]$ is written as
\begin{eqnarray}\label{matp^n}
\left(
 \begin{array}{cc}
 \chi_{p^n}(g) \psi_n^{-1}(g) &u_n(g)\\
 0& \psi_n(g)
 \end{array}
 \right).
\end{eqnarray}
Here $\chi_{p^n}$ denotes the mod $p^n$ cyclotomic character, $\psi_n$ is the mod $p^n$ reduction of $\psi$, and $u_n(g)=u(g)$ mod $p^n$.\\

 \begin{prop}\label{dims}
 Suppose that  $E$ has good ordinary reduction at $p$. Consider the following 4 cases.
\begin{itemize}
\item[$(A)$] $ a^j_p \not\equiv 1 \mod p$.
\item[$(B)$] $ a^j_p \equiv 1 \mod p$, $E$ has CM over an extension of $\mathbb{Q}_p$.
\item[$(C)$] $ a^j_p \equiv 1 \mod p$, $E$ does not have CM over an extension of $\mathbb{Q}_p$ and $\rho_{E, p}(G_{\mathbb{Q}_p})$ is not diagonalizable.
\item[$(D)$] $ a^j_p \equiv 1 \mod p$, $E$ does not have CM over an extension of $\mathbb{Q}_p$, $\rho_{E, p}(G_{\mathbb{Q}_p})$ is diagonalizable.
\end{itemize}
Then, in each case, the dimensions $\dim_Fp \frac{H^0(\mathbb{Q}_p, A^j_p)}{pH^0(\mathbb{Q}_p, A^j_p)}$ and $\mathrm{dim}_{\mathbb{Q}_p}H^0(\mathbb{Q}_p, V^j_p)$ are as in the table below.
\begin{table}[h]
  \centering
  \begin{tabular}{|c|c|c|} \hline
    \backslashbox{Case}{dimension} & $\dim_Fp \frac{H^0(\mathbb{Q}_p, A^j_p)}{pH^0(\mathbb{Q}_p, A^j_p)}$ & $\mathrm{dim}_{\mathbb{Q}_p}H^0(\mathbb{Q}_p, V^j_p)$ \\ \hline
    $(A)$ & 0 & 0  \\ \hline
    $(B)$ & 1 & 0  \\ \hline
    $(C)$ & 0 & 0  \\ \hline
    $(D)$ & 1 & 0  \\ \hline
\end{tabular}
\end{table}
 \end{prop}
(Proof of Proposition \ref{dims})

Here we also describe precise computations only for $\dim_Fp \frac{H^0(\mathbb{Q}_p, A^j_p)}{pH^0(\mathbb{Q}_p, A^j_p)}$. First we assume that $E$ has CM over some extension of $\mathbb{Q}_p$. For each positive integer $n$, we take $\overline{v_1}:= v_1$ mod $p^n$, $\overline{v_2}:=v_2$ mod $p^n$ as a basis of $E[p^n]$ over $\mathbb{Z}/p^n$, and we take a basis of $\Sym^j E[p^n]$ over $\mathbb{Z}/p^n$ as follows:
 \[
 w_0:=\overline{v_1}^{\otimes j}, w_1:=\overline{v_1}^{\otimes {j-1}} \otimes \overline{v_2}, \ldots , w_i:=\overline{v_1}^{\otimes {j-i}} \otimes \overline{v_2}^{\otimes i} , \ldots , w_j:=\overline{v_2}^{\otimes j}.
 \]
 Since $E$ has CM, we may assume that $u(g)$ in  (\ref{matp^n}) is 0 for all $g \in G_{\mathbb{Q}_p}$. Then we have 
 \[
 g(w_i)=g(\overline{v_1}^{\otimes {j-i}} \otimes \overline{v_2}^{\otimes i})=(\chi_{p^n}(g)\overline{v_1})^{\otimes {j-i}} \otimes \overline{v_2}^{\otimes i}=\chi_{p^n}(g)^{j-i} w_i\ \ (g \in G_{\mathbb{Q}_p^{\mathrm{ur}}}).
 \]
 For an element $x:=a_0w_0+a_1w_1+\cdots+a_jw_j\ (a_i \in \mathbb{Z}/p^n)$ in $\Sym^j E[p^n]$, the condition $x \in \Sym^j E[p^n]^{G_{\mathbb{Q}_p^{\mathrm{ur}}}}$ is equivalent to the condition  
 \[
 a_0\chi_{p^n}(g)^j=a_0,\ a_1\chi_{p^n}(g)^{j-1}=a_1, \ldots, a_{j-1}\chi_{p^n}(g)=a_{j-1},\ a_j \in \mathbb{Z}/p^n\ \ (\forall g \in G_{\mathbb{Q}_p^{\mathrm{ur}}}).
 \]
Since $j\leqslant p-2$, for each $i$ there exists $g \in  G_{\mathbb{Q}_p^{\mathrm{ur}}}$ such that $\chi_{p^n}(g)^{j-i}-1 \in (\mathbb{Z}/p^n)^{\times}$. Thus we have $\Sym^j E[p^n]^{G_{\mathbb{Q}_p^{\mathrm{ur}}}}=\mathbb{Z}/p^n \cdot w_j$ to get $\Sym^j E[p^n]^{G_{\mathbb{Q}_p}}=(\mathbb{Z}/p^n\cdot w_j)^{\mathrm{Frob}_p=1}$. We know $\mathrm{Frob}_p$ acts on $w_j$ via the character $\psi^j_n=\psi^j\mod p^n$ in (\ref{matp^n}). Since $\psi$ is an infinite order character, there exists a non-negative integer $s$ such that $\psi(\mathrm{Frob}_p)^j\equiv 1 \mod p^s$ and $\psi(\mathrm{Frob}_p)^j \not\equiv 1 \mod p^{s+1}$. Then 
\[
\Sym^j E[p^n]^{G_{\mathbb{Q}_p}} = \left\{\begin{array}{ll}
\mathbb{Z}/p^n w_j  & (n\leqslant s) \\  
p^{n-s}\mathbb{Z}/p^n w_j & (n\geqslant s+1).
\end{array}\right.
\]
 Especially if $s=0$, in other words if $a^j_p \equiv \psi(\mathrm{Frob}_p)^j \not\equiv 1 \mod p$, then we have $\Sym^j E[p^n]^{G_{\mathbb{Q}_p}}=0$\ for all $n$, and we get
 \[
(A_p^j)^{G_{\mathbb{Q}_p}} \cong \left\{\begin{array}{ll}
0  & (a^j_p \not\equiv 1 \mod p) \\  
\displaystyle\frac{1}{p^{s}}\mathbb{Z}/\mathbb{Z} & (1\leqslant s < \infty).\\
\end{array}\right.
\]
Thus we get the desired result in the case (B) and partially in the case (A) when $E$ has CM.

\vspace{1mm}
Next we consider the case where $E$ does not have CM over an extension of $\mathbb{Q}_p$. In this case, there exists a non-negative integer $m$ such that $\rho_{E, p^m}(G_{\mathbb{Q}_p})$ is diagonalizable and $\rho_{E, p^{m+1}}(G_{\mathbb{Q}_p})$ is not diagonalizable. For each $n$, we take the basis $w_0, w_1, \ldots , w_j$ for $\Sym^j E[p^n]$ over $\mathbb{Z}/p^n$ as in the previous argument. 

For any $n \leqslant m$, we may assume $u_n(g)=0$ in (\ref{matp^n}) for all $g \in G_{\mathbb{Q}_p^{\mathrm{ur}}}$. So we can imitate the argument in the case where $E$ has CM, and get $\Sym^j E[p^n]^{G_{\mathbb{Q}_p}}=(\mathbb{Z}/p^n\cdot w_j)^{\mathrm{Frob}_p=1}$. 

For $n\geqslant m+1$, we have $u_n(G_{\mathbb{Q}_p}) \neq 0$ and $u_n(G_{\mathbb{Q}_p}) \subset p^m\mathbb{Z}/p^n\mathbb{Z}$. We first consider $\Sym^j E[p^n]^{G_{\mathbb{Q}_p^{\mathrm{ab}}}}$. With respect to the basis $\{\overline{v_1}, \overline{v_2}\}$, the group $G_{\mathbb{Q}_p^{\mathrm{ab}}}$ acts on $E[p^n]$ via
\begin{eqnarray*}
\left(
 \begin{array}{cc}
 1 &u_n(g)\\
 0&1
 \end{array}
 \right).
\end{eqnarray*}
 For an element $x:=a_0w_0+a_1w_1+\cdots+a_jw_j\ (a_i \in \mathbb{Z}/p^n)$ in $\Sym^j E[p^n]$, we can show that the condition $x \in \Sym^j E[p^n]^{G_{\mathbb{Q}_p^{\mathrm{ab}}}}$ is equivalent to the condition
\begin{eqnarray}\label{nonCMcondi} 
a_0 \in \mathbb{Z}/p^n\mathbb{Z},\ a_1 u_n (g)=a_2 u_n (g)=\ldots=a_j u_n (g)=0\ \ ( \forall g\in G_{\mathbb{Q}_p^{\mathrm{ab}}})
\end{eqnarray}
by exactly the same computation to the one in the proof of Lemma \ref{comp} if we think $u_n$ as $\tau_{q, n}$. We have $u_n(G_{\mathbb{Q}_p^{\mathrm{ur}}})=p^m\mathbb{Z}/p^n\mathbb{Z}$ by the definition of the integer $m$ and \cite[Lemma 3.5]{LA}. Here we study the image of $G_{\mathbb{Q}_p^{\mathrm{ab}}}$ under the map $u_n$.
\begin{lem}\label{imab}
For $n \geqslant m+1$, $u_n(G_{\mathbb{Q}_p^{\mathrm{ab}}})=p^m\mathbb{Z}/p^n\mathbb{Z}$.
\end{lem}
(Proof of Lemma \ref{imab})

 When $n=m+1$, $u_{m+1}(G_{\mathbb{Q}_p^{\mathrm{ab}}})=p^m\mathbb{Z}/p^{m+1}\mathbb{Z}\ \text{or}\ 0$ since  $u_{m+1}(G_{\mathbb{Q}_p^{\mathrm{ab}}})$ forms an additive group. If $u_{m+1}(G_{\mathbb{Q}_p^{\mathrm{ab}}})=0$, then $\rho_{E, m+1}(G_{\mathbb{Q}_p^{\mathrm{ur}}})$ is abelian but we can show that this can not be happen using \cite[Lemma 3.5]{LA} and the definition of the integer $m$. Thus we have
$u_{m+1}(G_{\mathbb{Q}_p^{\mathrm{ab}}})=p^m\mathbb{Z}/p^{m+1}\mathbb{Z}$. For $n \geqslant m+1$, taking compatible bases of $E[p^n]$ for all $n$ as in \cite[Lemma 3.5]{LA}, we have $u_n(g) \equiv u_{m+1}(g) \pmod{p^{m+1}}$ for $g \in G_{\mathbb{Q}_p^{\mathrm{ab}}}$ and obtain $u_n(G_{\mathbb{Q}_p^{\mathrm{ab}}})=p^m\mathbb{Z}/p^n\mathbb{Z}(=u_n(G_{\mathbb{Q}_p^{\mathrm{ur}}}))$. \hfill$\square$

 Hence if $x \in \Sym^j E[p^n]^{G_{\mathbb{Q}_p^{\mathrm{ab}}}}$, we have $a_1, a_2, \ldots, a_j \in p^{n-m}\mathbb{Z}/p^n\mathbb{Z}$ from (\ref{nonCMcondi}) and 
 \begin{eqnarray}\label{unrademo}
 a_1 u_n (g)=a_2 u_n (g)=\ldots=a_j u_n (g)=0 \ \ \text{in}\  \mathbb{Z}/p^{n}\mathbb{Z}
\end{eqnarray}
  still for $g \in G_{\mathbb{Q}_p^{\mathrm{ur}}}$ since $u_n(G_{\mathbb{Q}_p^{\mathrm{ur}}})=p^m\mathbb{Z}/p^n\mathbb{Z}$. 
  
  We next consider a condition on $a_0, a_1, \ldots, a_j$ such that $x=a_0 w_0+\cdots +a_j w_j\in \Sym^j E[p^n]^{G_{\mathbb{Q}_p^{\mathrm{ur}}}}$. For $g \in G_{\mathbb{Q}_p^{\mathrm{ur}}}$ and $i$ with $1\leqslant i\leqslant j$,
\begin{eqnarray*}
g(a_i w_i)=a_i g(\overline{v_1}^{\otimes{j-i}} \otimes \overline{v_2}^{\otimes{i}})& = &a_i (\chi_{p^n}(g) \overline{v_1})^{\otimes{j-i}}\otimes (u_n(g)\overline{v_1}+\overline{v_2})^{\otimes{i}}\\
& = & a_i\chi_{p^n}^{j-i}(g)\overline{v_1}^{\otimes{j-i}} \otimes\left( \sum_{k=0}^{i} \binom{i}{k} u_n(g)^{i-k} \overline{v_1}^{\otimes{i-k}} \otimes \overline{v_2}^{\otimes k}\right)\\
&= & a_i \chi_{p^n}^{j-i}(g) \overline{v_1}^{\otimes{j-i}} \otimes \overline{v_2}^{\otimes i} =a_i \chi_{p^n}^{j-i}(g) w_i.
\end{eqnarray*}
Here we use (\ref{unrademo}) in the fourth equality for $i$ with $1\leqslant i\leqslant j$ since $x \in \Sym^j E[p^n]^{G_{\mathbb{Q}_p^{\mathrm{ur}}}} \subset \Sym^j E[p^n]^{G_{\mathbb{Q}_p^{\mathrm{ab}}}}$. For $i=0$, we also have $g(a_0 w_0)=a_0 \chi_{p^n}^{j}(g) w_0$ for $g \in G_{\mathbb{Q}_p^{\mathrm{ur}}}$. So the condition $x \in \Sym^j E[p^n]^{G_{\mathbb{Q}_p^{\mathrm{ur}}}}$ is equivalent to the condition
 \begin{eqnarray}\label{condfornonCMpart}
 a_0\chi_{p^n}(g)^j=a_0, \ldots, a_{j-1}\chi_{p^n}(g)=a_{j-1},\ a_j \in p^{n-m}\mathbb{Z}/p^n\mathbb{Z}
 \end{eqnarray}
for all $g \in G_{\mathbb{Q}_p^{\mathrm{ur}}}$. Again we can take $g \in G_{\mathbb{Q}_p^{\mathrm{ur}}}$ such that $\chi_{p^n}(g)^{j-i}-1 \in (\mathbb{Z}/p^n)^{\times}$ for each $i$ to get $a_0=a_1=\ldots=a_{j-1}=0$. Thus we get $\Sym^j E[p^n]^{G_{\mathbb{Q}_p}}=(p^{n-m}\mathbb{Z}/p^n\cdot w_j)^{\mathrm{Frob}_p=1}$ for $n\geqslant m+1$. 

If $\rho_{E, p}(G_{\mathbb{Q}_p})$ is not diagonalizable, in other words if $m=0$, $\Sym^j E[p^n]^{G_{\mathbb{Q}_p}}=\Sym^j E[p^n]^{G_{\mathbb{Q}^{\mathrm{ur}}_p}}=0$ for all $n$ from the above computations and $({A_p^j})^{G_{\mathbb{Q}_p}}=0$. Thus we get the desired result in the case $(C)$ and partially in the case (A).
 
 If $m\geqslant 1$, from the above argument, we have 
\[
\Sym^j E[p^n]^{G_{\mathbb{Q}_p}} = \left\{\begin{array}{ll}
(\mathbb{Z}/p^n w_j)^{\mathrm{Frob}_p=1}  & (n\leqslant m) \\  
(p^{n-m}\mathbb{Z}/p^n\mathbb{Z} \cdot w_j)^{\mathrm{Frob}_p=1} & (n\geqslant m+1).
\end{array}\right.
\]
We know that $\mathrm{Frob}_p$ acts on $\Sym^j E[p^n]^{G_{\mathbb{Q}_p^{\mathrm{ur}}}}$ via the character $\psi^j_n=\psi^j\mod p^n$ for all $n$. We again take a non-negative integer $s$ such that $\psi(\mathrm{Frob}_p)^j \equiv 1 \mod p^s$ and $\psi(\mathrm{Frob}_p)^j  \not\equiv 1 \mod p^{s+1}$.
If $s=0$, in other words if $a^j_p \not\equiv 1\mod p$, then we have $\Sym^j E[p^n]^{G_{\mathbb{Q}_p}}=0$ for all $n$. If $s >0$ we have 
\[
\Sym^j E[p^n]^{G_{\mathbb{Q}_p}} \cong \left\{\begin{array}{ll}  
 \mathbb{Z}/p^n w_j & (n\leqslant \mathrm{min}\{m, s\} )\\
p^{n- \mathrm{min}\{m, s\}}\mathbb{Z}/p^n w_j & (n > \mathrm{min}\{m, s\}).
\end{array}\right.
\]
Thus we get
\[
(A_p^j)^{G_{\mathbb{Q}_p}} \cong \left\{\begin{array}{ll}
0  & (a^j_p \not\equiv 1 \mod p) \\  
\displaystyle\frac{1}{p^{\mathrm{min}\{m, s\}}}\mathbb{Z}/\mathbb{Z} & (1\leqslant s < \infty).
\end{array}\right.
\]
So finally, we have the desired result in the case $(D)$ and the case $(A)$ completely.

\hfill$\square$
\section{Non-injectivity of the restriction map $\mathrm{Res}_{\mathbb{Q}_p^{\mathrm{ur}}/\mathbb{Q}_p}$}
From (\ref{ineq}) and Proposition \ref{dims}, we deduce Proposition \ref{dimofResp} in the case $(A)$ and $(C)$, and the main theorem follows. Since we assume $(b^{\prime})$ in Theorem \ref{main}, $(D)$ in the table in Proposition \ref{dims} does not occur. For the case $(B)$, we prove the following proposition.
\begin{prop}\label{noninj}
 In the case $(B)$, the restriction map 
 \[
 \mathrm{Res}_{\mathbb{Q}_p^{\mathrm{ur}}/\mathbb{Q}_p}:\iota^{-1}(H_f^1(\mathbb{Q}_{p}, A^j_p)[p]) \rightarrow H^1(\mathbb{Q}_p^{\mathrm{ur}}, \Sym^j E[p])
 \]
  is not injective.
 \end{prop}
 From this proposition, also in the case $(B)$, we can deduce Proposition \ref{dimofResp} and the main theorem follows. Thus the main theorem follows in all possible cases $(A), (B)$ and  $(C)$ under the assumptions in Theorem \ref{main}.\\

\noindent
(Proof of Proposition \ref{noninj})\\
We have the  following commutative diagram
\[
  \xymatrix{
 0 \ar[r]^{}&   \displaystyle  \frac{H^0(\mathbb{Q}_p, A^j_p)}{pH^0(\mathbb{Q}_p, A^j_p)} \ar[d]\ar[r]  & \iota^{-1}(H^1_f(\mathbb{Q}_{p}, A^j_p)[p]) \ar[d]^{\mathrm{Res}_{\mathbb{Q}_p^{\mathrm{ur}}/\mathbb{Q}_p}} \\
 0 \ar[r]^{}&   \displaystyle  \frac{H^0(\mathbb{Q}^{\mathrm{ur}}_p, A^j_p)}{pH^0(\mathbb{Q}^{\mathrm{ur}}_p, A^j_p)} \ar[r] &  H^1(\mathbb{Q}^{\mathrm{ur}}_p, \Sym^j E[p]) . 
  }
\]
From the table in Proposition \ref{dims}, we have $\dim_Fp\left(\frac{H^0(\mathbb{Q}_p, A^j_p)}{pH^0(\mathbb{Q}_p, A^j_p)}\right)=1$ in the case $(B)$. Since we have already computed $\Sym^j E[p^{\infty}]^{G_{\mathbb{Q}_p^{\mathrm{ur}}}}$ in the case in the proof of Proposition \ref{dims}, we also get $\dim_Fp\left(\frac{H^0(\mathbb{Q}^{\mathrm{ur}}_p, A^j_p)}{pH^0(\mathbb{Q}^{\mathrm{ur}}_p, A^j_p)}\right)=0$. So by the above commutative diagram, the dimension of the kernel of $\mathrm{Res}_{\mathbb{Q}_p^{\mathrm{ur}}/\mathbb{Q}_p}$ is at least $1$. Thus the proposition follows.

\hfill$\square$

\subsection*{acknowledgement}
The author would like to thank his supervisor Professor Masato Kurihara heartily for guiding him to the topics in this paper, careful reading of the draft of this paper and helpful discussions with him. He also would like to express his sincere gratitude to Professor Dipendra Prasad and Professor Sudhanshu Shekhar for many valuable comments on his work. Thanks are also due to Ryotaro Sakamoto for answering his questions precisely concerning Proposition 2. Lastly, he is also grateful to the referees for their careful reading of his draft. This research was supported by JSPS KAKENHI Grant Number 21J13502.

\end{document}